\documentclass[a4paper,12pt]{article}

\usepackage{amssymb}
\usepackage{epsfig}
\usepackage{amsfonts}
\usepackage{amsmath}
\usepackage{euscript}
\usepackage{amscd}
\usepackage{amsthm}
\DeclareMathAlphabet{\mathpzc}{OT1}{pzc}{m}{it}

\newtheorem{thm}{Theorem}[section]

\newcommand{\m}{\mathpzc{m}}
\newcommand{\p}{\mathpzc{p}}

\newcommand{\bZ}{\mathbb Z}
\newcommand{\bQ}{\mathbb Q}
\newcommand{\bR}{\mathbb R}
\newcommand{\bC}{\mathbb C}

\newcommand{\A}{\mathbb A}

\newcommand{\Ga}{\mathbb G}

\newcommand{\dk}{\operatorname{DK}}

\title{The Zariski Cancellation Problem and related problems in Affine Algebraic Geometry}
\author{Neena Gupta\\
{\small{\it Statistics and Mathematics  Unit, Indian Statistical Institute,}}\\
{\small{\it 203 B.T. Road, Kolkata 700 108, India}}\\
{\small{\it e-mail : neenag@isical.ac.in, rnanina@gmail.com}}\\
}

\begin{document}

\date{}
\maketitle

\begin{abstract}
In this article, we shall discuss the solution to the Zariski Cancellation Problem
in positive characteristic, various approaches taken so far towards
the possible solution in characteristic zero, and several other questions
related to this problem.
\smallskip

\end{abstract}

\section{Introduction}

\begin{quote}
``Polynomials and power series \\
\hspace*{2mm}May they forever rule the world.''
\hfill---Shreeram S. Abhyankar, 1970 (\cite{AAG})
\end{quote} 
Right from the beginning of the nineteenth century, mathematicians have been involved in studying polynomial rings (over $\bC$ and over $\bR$).  
Some of the early breakthroughs on polynomial rings have led to the foundation of the subject Commutative Algebra. 
One such result is the Hilbert Basis Theorem, a  landmark result on the finite generation of ideals, which solved a 
central problem on invariant theory. This was followed by the Hilbert Nullstellensatz which connects affine varieties 
(zero locus of a set of polynomials) with rings of regular functions on varieties 
and thus enables one to make use of the algebraic machinery of commutative  algebra to study geometric properties of varieties.

Affine Algebraic Geometry deals with the study of affine spaces (and certain closed subspaces), equivalently, polynomial rings (and certain quotients). 
There are many fundamental problems on polynomial rings  which can be formulated in an elementary mathematical 
language but whose solutions remain elusive.  Any significant progress requires development of new and powerful methods and their ingenious applications. 

One of the most challenging problems in Affine Algebraic Geometry is the  
Zariski Cancellation Problem  (ZCP)  on polynomial rings (Question 1$^\prime$ below).
In this article, we shall discuss the solution to the ZCP in positive characteristic,  various approaches taken so far  towards  
the possible solution in characteristic zero, and several other questions related to this problem.  
For a survey on problems in Affine Algebraic Geometry one may  look at \cite{Kr}, \cite{M} and \cite{R1}.

Throughout the article,  all rings will be assumed to be commutative with unity and $k$ will denote a field. 
For a ring $R$, $R^*$ will denote the group of units of $R$.
We shall use the notation $R^{[n]}$  for a polynomial ring in $n$ variables over a commutative ring $R$.
Thus, $E= R^{[n]}$ will mean that $E= R[t_1, \ldots, t_n]$ for some elements  
$t_1, \ldots, t_n$ in $E$ which are algebraically independent over $R$.
Unless otherwise stated, capital letters like $X_1, X_2,\dots, X_n,Y_1, \dots, Y_m, X, Y, Z,T$
will be used as variables of polynomial rings.

\section{Cancellation Problem}

Let $A$ be an affine (finitely generated) algebra over a field $k$. 
The $k$-algebra A is said to be  cancellative (over $k$) if, for any $k$-algebra $B$,
$A[X] \cong_k B[X]$ implies that $A\cong_k B$.  A natural question in this regard is:  
which affine domains are cancellative? More precisely:

\medskip

\noindent
{\bf Question 1}. Let $A$ be an affine algebra over a field $k$. 
Suppose that $B$ is a $k$-algebra such that the polynomial 
rings $A[X]$ and $B[X]$ are isomorphic as $k$-algebras.
Does it follow that $A\cong_k B$? In other words, 
is the $k$-algebra $A$ cancellative? 

\medskip

A special case of Question 1, famously known as the Zariski Cancellation Problem,  
asks whether affine spaces are cancellative, i.e.,  whether any polynomial ring in $n$ variables over a field $k$ is cancellative.  More precisely:

\medskip

\noindent
{\bf Question 1$^\prime$}. Suppose that $B$ is an affine $k$ algebra satisfying 
$B[X] \cong_k k[X_1, \dots, X_{n+1}]$ for some positive integer $n$. 
Does it follow that $B\cong_k k[X_1, \dots, X_n]$?
In other words, is the polynomial ring $k[X_1, \dots, X_{n}]$ cancellative?

\medskip

S.S. Abhyankar, P. Eakin and W.J. Heinzer have shown that any domain $A$ 
of transcendence degree one over any field $k$ is cancellative (\cite{AEH}). 
 In fact they showed that, for any UFD $R$, the polynomial ring $R[X]$ is cancellative over $R$. 
 This was further generalised by E. Hamann to a ring $R$  which either contains $\bQ$ or is a seminormal domain (\cite{H}).  

In 1972,  M. Hochster demonstrated  the first counterexample to Question 1 (\cite{Ho}). 
His example,  a four dimensional ring over the field of real numbers $\bR$, 
 is based on the fact that the projective module defined by the tangent bundle over the real sphere with 
 coordinate ring $S= \bR[X,Y,Z]/(X^2+Y^2+Z^2-1)$ is stably free but not a free $S$-module. 
 
 One of the major breakthroughs in 1970's was the establishment of 
 an affirmative answer to Question 1$^\prime$ for the case $n=2$. 
 This was proved over a field of characteristic zero  by T. Fujita, M. Miyanishi and T. Sugie (\cite{F}, \cite{MS}) and
 over perfect fields of arbitrary characteristic by P. Russell  (\cite{R}). 
 Later, it has been shown that even the hypothesis of perfect field can be dropped (\cite{BG3}).
A simplified proof of the cancellation property of $k[X, Y]$ 
 for an algebraically closed field $k$ is given by A. Crachiola and L. Makar-Limanov in \cite{CML}.
 
Around 1989, W. Danielewski  (\cite{Da}) constructed explicit two dimensional affine domains 
over the field of complex numbers $\bC$ which are not cancellative over $\bC$. 
New examples of non-cancellative surfaces over any field $k$ have been studied in  \cite{NS1}.
This addresses the cancellation problem, as formulated in Question 1, for all dimensions. 

In \cite{G1} and \cite{G3}, the author  settled the Zariski Cancellation Problem (Question 1$^\prime$) 
completely for affine spaces in positive characteristic.
She has first shown in \cite{G1}  that a certain threefold constructed by 
T. Asanuma is a counterexample to the ZCP in positive characteristic
for the affine three space. Later in \cite{G2}, she  studied a general threefold of the
form $x^my= F(x,z,t)$, which includes the Asanuma threefold as well as the 
famous Russell cubic defined below.
A major theorem of \cite{G2} is stated as Theorem \ref{ngc} of this article. 
In \cite{G3}, using a modification of the theory developed in \cite{G2},
she constructed a family of examples which are counterexamples to the ZCP
in positive characteristic in all dimensions greater  than $2$. 
The ZCP is still a challenging problem in characteristic zero. 
A few candidate counterexamples are discussed below.

\medskip

\noindent
{\bf The Russell cubic}: 

\medskip

Let $A= \bC[X,Y,Z,T]/(X^2Y+X+Z^2+T^3)$, $V= {\rm Spec}~A$ and let $x$ denote the image of $X$ in $A$. 
The ring $A$,  known as the Russell cubic, is one of the simplest examples of the Koras-Russell threefolds,
a family of threefolds which arose in the context of the problem
of determining whether there exist non-linearisable ${\bC}^{*}$-actions on $\bC^3$. 
It was an exciting open problem for some time whether $A\cong \bC^{[3]}$. 
It was first observed that the ring $A$ (respectively the variety $V$) has several properties 
in common with $\bC^{[3]}$ (respectively $\bC^3$), for instance, 

\medskip

\noindent
(i)  $A$ is a regular UFD.\\
(ii) There exists  an injective $\bC$-algebra homomorphism from 
$A$ to $\bC^{[3]}$. Note that  $\bC^{[3]} \hookrightarrow A$.\\
(iii) The variety $V$ is homeomorphic (in fact diffeomorphic) to $\bR^6$.\\
(iv) $V$ has logarithmic Kodaira dimension $-\infty$. 

\medskip
 
These properties appeared to provide evidence in favour of the surmise that $A \cong \bC^{[3]}$.
The establishment of an isomorphism between $A$ and $\bC^{[3]}$ would have led to counterexamples
to the ``Linearisation Conjecture'' on  $\bC^3$ (stated in \cite{Kam}) and the Abhyankar-Sathaye Conjecture for $n=3$
(stated in Section 5 of the present article).
For, if $A$ were isomorphic to $\bC^{[3]}$, as was then suspected, it would have  
shown the existence of non-linearisable $\bC^*$-actions on $\bC^3$. 
Moreover note that 

\medskip

\noindent
(v) $A/(x-\lambda) = \bC^{[2]}$ for every $\lambda \in \bC^*$.  \\
(vi) $A/(x) \neq \bC^{[2]}$.

\medskip

Therefore, if $A$ were isomorphic to $\bC^{[3]}$, then 
property (vi) would show that $x-\lambda$  cannot be a coordinate in $A$ for any $\lambda $ and then, by property (v), 
it  would have yielded a counterexample to the
Abhyankar-Sathaye Conjecture for $n=3$.

However, L. Makar-Limanov proved (\cite{MLR}) that $A \neq \bC^{[3]}$;
for this result he introduced  a new invariant which distinguished between $A$ and $\bC^{[3]}$.
This invariant, which he had named AK-invariant, is now named Makar-Limanov invariant
and is denoted by $ML$. It  is defined in Section 3.
Makar-Limanov proved that

\medskip

\noindent
(vii) $ML(A)= \bC[x]$ (Makar-Limanov \cite{MLR}).
\medskip

However, the Makar-Limanov invariant of $\bC^{[n]}$ is $\bC$ for any integer $n \ge 1$.
Thus $A\ncong \bC^{[3]}$. Subsequently, other Koras-Russell threefolds were shown to be not isomorphic to the polynomial ring. 
Eventually Kaliman--Koras--Makar-Limanov--Russell proved that every $\bC^*$-action on
$\bC^{3}$ is linearisable (cf. \cite{KKMLR}). 

Now for ZCP in characteristic zero, a crucial question, still open, is whether $A^{[1]}=\bC^{[4]}$.
Because if $A^{[1]}=\bC^{[4]}$, then $A$ would be a counterexample to the ZCP in characteristic zero for $n=3$.
In this context, the following results have been proved:

\medskip

\noindent
 (viii) $ML(A^{[1]})=\bC$ (Dubouloz \cite{Dub}). \\
 (ix) $V$ is $\A^1$-contractible (Dubouloz-Fasel \cite{DuFa}, also see \cite{HKO}, \cite{PDO}). 

\medskip

Note that $A^{[1]}=\bC^{[4]}$ would imply  that $ML(A^{[1]})= \bC$ and  A. Dubouloz's result (viii) shows that the latter indeed holds.
On the other hand, A. Asok had suggested a program for showing that 
the variety $V$ is not $\A^1$-contractible and hence A is not a stably polynomial ring (see \cite{HKO}). 
However, M. Hoyois, A. Krishna and P.A. Østvær have proved (\cite{HKO}) that a step in his program does not hold for $V$. 
They had further shown that $V$ is stably $\A^1$-contractible. 
In a remarkable paper (\cite{DuFa}), A. Dubouloz and J. Fasel  have established that $V$ is in fact $\A^1$-contractible
  which seems to provide further evidence in favour of $A^{[1]}=\bC^{[4]}$.
The variety $V$ is in fact the first example of an $\A^1$-contractible threefold which is not algebraically isomorphic to $\bC^3$.

\medskip

\medskip

\noindent
{\bf Non-rectifiable epimorphisms and Asanuma's rings}:

 Let $m\le n$ be two integers.
A $k$-algebra epimorphism $\phi: k[X_1, \dots, X_n] \twoheadrightarrow k[Y_1, \dots, Y_m]$ is said to be {\it rectifiable}
if there exists a $k$-algebra automorphism $\psi$ of $k[X_1, \dots, X_n]$ 
such that ${{\phi}}_{\circ} \psi (X_i)= Y_i$ for $1\le i\le m$ and ${{\phi}}_{\circ} \psi(X_j)= 0$
for $m+1\le j\le n$. Equivalently, over an algebraically closed field $k$,  a $k$-embedding 
$\Phi : \A^m_k \hookrightarrow \A^n_k$ is said to be {\it rectifiable} if there 
exists an automorphism ${\Psi}$ of $\A^n_k$ such that $\Psi_\circ{{\Phi}}$
is the canonical embedding mapping
$(y_1, \dots, y_m) \to (y_1, \dots, y_m, 0,\dots, 0)$. 

A famous theorem of S.S. Abhyankar-T. Moh and M. Suzuki proves that 
 any epimorphism $\phi: k[X,Y] \to k[T]$ is rectifiable  in characteristic zero (\cite{AM},  \cite{Suz}). 
On the other hand, in positive characteristic,  there exist non-rectifiable epimorphisms from $ k[X,Y]$ to $k[T]$ 
(B. Segre  \cite{Se}, M. Nagata  \cite{Na}).
It is an open problem whether there exist non-rectifiable epimorphisms over the field of complex numbers (See \cite{DuG}).

T. Asanuma has described an explicit method for constructing affine rings
which are stably polynomial rings, by making use of non-rectifiable epimorphisms (\cite{A3}, also see \cite[Proposition 3.7]{DuG}).  
Such rings are considered to be  potential  candidates for counterexamples to the ZCP.
For instance,  when $k$ is of positive characteristic,  non-rectifiable epimorphisms from $k[X,Y]$ to $k[T]$ yield counterexamples to the ZCP.

Let $\phi: \bR[X,Y,Z] \to \bR[T]$ be defined by 
$$
\phi(X)=   T^3-3T, ~~\phi(Y)=T^4-4T^2,~~\phi(Z)=T^5-10T.
$$
A.R. Shastri constructed the above epimorphism $\phi$ and proved
 that it defines a non-rectifiable (polynomial) embedding 
of the trefoil knot  in $\A^3_{\bR}$ (\cite{Sh}).  
Using a result of J.P. Serre (\cite[Theorem 1, p. 281]{La}) one knows that ker$(\phi)= (f,g)$ for some $f, g \in k[X,Y,Z]$.
Using $f$ and $g$, Asanuma constructed the ring  $B = \bR[T][X,Y,Z,U,V]/(T^d U-f, T^dV-g)$ and proved that
 $B^{[1]}= \bR[T]^{[4]}=\bR^{[5]}$ (cf. \cite[Corollary 4.2]{A3}). 
He asked (\cite[Remark 7.8]{A3}):

\medskip

\noindent
{\bf Question 2.} Is $B= \bR^{[4]}$?

\medskip

The interesting aspect of the question is that once the problem gets solved then, 
 irrespective of whether the  answer is `Yes' or `No', 
that is, either way, one would have solved a major problem in Affine Algebraic Geometry. For:

\smallskip

\noindent
If $B=  \bR^{[4]}$, then there exist non-linearizable $\bR^*$-actions on the affine four space  $\A^{4}_{\bR}$. 

\smallskip

\noindent
If $B \neq  \bR^{[4]}$, then clearly $B$ is a counter-example to the ZCP!!

\section{Characterisation Problem}

The Characterisation Problem  in affine algebraic geometry seeks  
 a ``useful characterisation'' of 
the polynomial ring or, equivalently (when the ground field is algebraically closed)
an affine $n$-space.
 For instance, the following two results give respectively an algebraic and a topological  characterisation of $k^{[1]}$ (or
 $\A^{1}_{\bC}$).

  \begin{thm}
Let $k$ be an algebraically closed field of characteristic zero.
Then  the polynomial ring $k^{[1]}$ is the only one-dimensional affine UFD with $A^*=k^*$.
\end{thm}

\begin{thm}
Let $k$ be the field of complex numbers $\bC$.  Then the affine line $\A^{1}_{\bC}$ 
is the only acyclic normal curve. 
\end{thm}

While the Characterisation Problem is one of the most important problems in affine algebraic geometry in its own right, it 
is also closely related to some of the challenging open problems on the affine space like the ``Cancellation Problem''. 
For instance, each of the above characterisations of $k^{[1]}$ immediately solves the Cancellation Problem 
in dimension one: $A^{[1]} =k^{[2]} \implies A= k^{[1]}$.
The complexity of the characterisation problem increases with the dimension of the rings.

In his attempt to solve the Cancellation Problem for the affine plane, C.P. Ramanujam obtained a remarkable topological
characterisation of the affine plane $\bC^2$ (1971, \cite{Ram}). He proved that
\begin{thm}
$\bC^2$ is the only contractible smooth surface which is simply connected at infinity. 
\end{thm}

Ramanujam also constructed contractible surfaces which are not isomorphic to $\bC^2$. 
Soon, M. Miyanishi (1975, \cite{Miya}) obtained an algebraic characterization of the polynomial ring
$k^{[2]}$. He proved that
\begin{thm}
Let $k$ be an algebraically closed field of characteristic zero and $A$ be a two dimensional affine factorial domain over $k$.
Then $A=k^{[2]}$ if and only if it satisfies the following:
\begin{enumerate}
\item[\rm(i)] $A^*=k^*$.
\item[\rm(ii)] There exists an element $f \in A$ and a subring $B$ of $A$ such that $A[f^{-1}]={B[f^{-1}]}^{[1]}$. 
\end{enumerate}
\end{thm}

This algebraic characterisation was used by  T. Fujita, M. Miyanishi and T. Sugie (1979-\cite{F}, 1980-\cite{MS})
to solve the Cancellation Problem for $k[X,Y]$. 
In 2002 (\cite{Gur}), using methods of Mumford and Ramanujam, R.V. Gurjar  gave a 
topological proof of the cancellation property of $\bC[X,Y]$. 


Remarkable characterisations of the affine three space were obtained by Miyanishi (1984, \cite{MN2}) and 
Kaliman (2002, \cite{Ka};  also see \cite{M} for a beautiful survey).  We state below the version of Kaliman. 

\begin{thm}\label{c3}
Let $A$ be a three dimensional  smooth factorial affine domain over the field of 
complex numbers $\bC$. Let $X= {\rm Spec}~A$. Then $A= \bC^{[3]}$ if and only if it satisfies the following:
\begin{enumerate}
\item[\rm(i)] $A^*=\bC^*$.
\item[\rm(ii)] $H_3(X, \bZ)=0$, or $X$ is contractible. 
\item[\rm(iii)] $X$ contains a cylinderlike open set $V$ such that $V \cong U \times \A^2$ for some curve $U$ and each irreducible component of the complement $X\setminus V$ has at most isolated singularities.   
\end{enumerate}
\end{thm}

When $A^{[1]}= \bC^{[4]}$,  it is easy to see that $A$ possesses the properties (i) and (ii) of Theorem \ref{c3}.
Thus, by Theorem \ref{c3}, the ZCP  for $\bC^{[3]}$ reduces to examining whether the condition (iii) 
necessarily holds for a $\bC$-algebra $A$ satisfying $A^{[1]}= \bC^{[4]}$.

In \cite{NDG}, we have obtained another characterisation of the affine three space using 
certain invariants of an affine domain defined by locally nilpotent derivations.  We state it below.

\medskip

\noindent
{\bf Locally nilpotent derivations and a characterisation of $\bC^{[3]}$}

\medskip

Let $B$ be an affine domain over a field $k$ of characteristic zero.  
A $k$-linear derivation $D$ on $B$ is said to be a locally nilpotent derivation if, 
for any $a \in B$ there exists an integer $n$ (depending on $a$) satisfying $D^n (a)=0$. 
 Let $LND(B)$ denote the set of all locally nilpotent 
$k$-derivations of $B$ and let  
$$
LND^{*}(B)=\{D \in LND(B)~|~Ds=1~ \text{for some} ~ s \in B\}.
$$ 
Then we define 
$$
ML(B):=\bigcap_{D \in LND(B)} Ker~D ~~\text{and}~~ML^{*}(B):= \bigcap_{D \in LND^{*}(B)} Ker~D.
$$ 
$ML(B)$, introduced by Makar-Limanov. is now called the Makar-Limanov invariant of $B$. 
$ML^*(B)$  is introduced by G. Freudenburg  in \cite[p. 237]{Fr}.
We call it the Makar-Limanov--Freudenburg invariant or ML-F invariant.
If $LND^{*}(B)=\emptyset$, we define 
$ML^{*}(B)$ to be $B$. We have obtained the following theorem \cite[Theorem 4.6]{NDG}.

\begin{thm}
Let $A$ be a three dimensional affine factorial domain over an algebraically closed field $k$ of characteristic zero. 
Then the following are equivalent:
\begin{enumerate}
 \item [\rm (I)]$A=k^{[3]}$.
 \item [\rm (II)]$ML^{*}(A)=k$.
 \item [\rm (III)]$ML(A)=k$ and $ML^{*}(A) \neq A$.
\end{enumerate}
\end{thm} 

A similar result has also been proved in dimension two under weaker hypotheses  (\cite[Theorem 3.8]{NDG}).
 The above characterisation of the affine three space does not extend to higher dimensions (\cite[Example 5.6]{NDG}.
So far, no suitable characterisation of the affine $n$-space for $n \ge 4$ is known to the author.  

\section{Affine Fibrations}\label{fibration}

Let $R$ be a commutative ring. A fundamental theorem of Bass-Connell-Wright and Suslin (\cite{BCW}, \cite{Su}) 
on the structure of locally polynomial algebras states that:

\begin{thm}\label{bcw}
Let $A$ be a finitely presented algebra over a ring $R$. Suppose that for each maximal 
ideal $\m$ of $R$, $A_{\m}= R_{\m}^{[n]}$ for some integer $n \ge 0$. 
Then $A \cong Sym_R{(P)}$ for some finitely generated projective $R$-module $P$ of rank $n$. 
\end{thm}

Now for a prime ideal $P$ of $R$, let $k(P)$ denote the residue field $R_P/ P R_P$.
The area of affine fibrations seeks to derive information about 
the structure and properties of an $R$-algebra $A$ from the  
information about the fibre rings $A \otimes_R k(P) (= A_P/P A_P)$ of $A$ 
at the points $P$ of the prime spectrum of $R$, i.e.,
at the prime ideals $P$ of $R$.

{\it An $R$-algebra $A$ is said to be an $\A^n$-fibration over $R$
if $A$ is a finitely generated flat $R$-algebra and 
for each prime ideal $P$ of $R$, $A \otimes_R k(P)= k(P)^{[n]}$.}

The most important problem on $\A^n$-fibrations,
due to B. Ve\v{\i}sfe\v{\i}ler and I.V. Dolga\v{c}ev \cite{DW},
can be formulated as follows:

\medskip

\noindent
{\bf Question 3.} Let $R$ be a Noetherian domain of dimension $d$ and $A$ be 
an $\A^n$-fibration over~$R$. 


(i)  If $R$ is regular, is $A\cong {\rm Sym}_R (Q)$ for some  projective module $Q$ over $R$?
(In particular, if $R$ is regular local, is then $A=R^{[n]}$?) 

(ii) In general, what can one say about the structure of $A$?


\medskip

Question 3 is considered a  hard problem.
When $n=1$, it has an affirmative answer for all $d$.
This has been established in the works of T. Kambayashi, 
M. Miyanishi and David Wright (1978, \cite{KM} and 1985, \cite{KW}).
Their results were further refined by A.K. Dutta 
who showed  that it is enough to assume the fibre conditions only on generic and co-dimension one fibres (\cite{D}; also see \cite{BD1}, \cite{DuO}, \cite{BDO}). 

In case $n=2$, $d=1$ and $R$ contains the field of rational numbers, 
an important theorem of A. Sathaye (1983, \cite{S}) 
gives an affirmative answer to Question 3(i). 
To prove this theorem Sathaye first generalised the Abhyankar-Moh expansion techniques originally developed over $k[[x]]$ to 
$k[[x_1, \dots, x_n]]$ (\cite{Sai}). The expansion techniques were used by Abhyankar-Moh to prove their famous Epimorphism theorem.
The generalised expansion techniques were further developed by Sathaye (\cite{S3}) to 
prove a conjecture of D. Daigle and G. Freudenburg. 
The result was a crucial step in Daigle-Freudenburg's theorem 
that the kernel of any triangular derivation
of $k[X_1, X_2, X_3, X_4]$ is a finitely generated $k$-algebra (\cite{DF1}).

When the residue field of $R$ is of positive characteristic,
T.~Asanuma has shown in \cite[Theorem 5.1]{A} 
that Question 3(i) has negative answer
for $n=2$, $d=1$ and 
the author has generalised Asanuma's ring (\cite{G3}) 
to give a negative answer to Question 3(i) 
for $n=2$ and any $d>1$  (also see \cite{G5}).
In Theorem \ref{ngc}, the author proved that in a special situation 
$\A^2$-fibration is indeed trivial.

However, if $n=2$, $d=2$ and $R$ contains the field of rational numbers,  Question 3(i) is an open problem.
A candidate counterexample is discussed in Section \ref{bdex}.

In the context of Question 3(ii), a deep work of 
T. Asanuma (1987, \cite{A}) 
provides a stable structure theorem for $A$.
As a consequence of Asanuma's structure theorem, it follows that 
if $R$ is regular local, then there exists an integer $m\ge 0$ such that
$A^{[m]}= R^{[m+n]}$. Thus it is very tempting to look for possible 
counterexamples to the affine fibration problem in order to obtain 
possible counterexamples  to the ZCP in characteristic zero.  
One can see \cite{BDi}, \cite{Du2}, \cite{Du3}, \cite{DF} and \cite[Section 3.1]{DuG} for more results on affine fibrations. 

So far we have considered affine fibrations where the fibre rings are polynomial rings. 
Bhatwadekar and Dutta have obtained some nice results on rings whose fibre rings are of the form $k[X, 1/X]$ (\cite{BDA1}, \cite{BDA2}).
Later Bhatwadekar, the author and A.M.  Abhyankar studied rings whose fibre rings are Laurent polynomial algebras or rings of the form 
$k[X, 1/f(X)]$, or of the form $k[X,Y, 1/(aX+b), 1/(cY+d)]$ for some $a, b, c, d\in k$ (\cite{BG1}, \cite{BG2}, \cite{AsB1}, \cite{AsB2}, \cite{G}).  
One of the results provides a Laurent polynomial analogue of Theorem \ref{bcw} and affine fibration problem Question 3.
More generally \cite[Theorems A and C]{BG2}:

\begin{thm}
Let $R$ be a Noetherian normal domain with field of fractions $K$ and $A$ be a faithfully flat 
$R$-algebra such that
\begin{enumerate}
\item[\rm(i)] $A\otimes_RK\cong K[X_1, \frac{1}{X_1}, \dots, X_n, \frac{1}{X_n}]$ 
\item[\rm(ii)] For each height one prime ideal $P$ of $R$, $A \otimes_R k(P) \cong k(P)[X_1, \frac{1}{X_1}, \dots, X_n, \frac{1}{X_n}]$ 
\end{enumerate}
Then $A$ is a locally Laurent polynomial algebra in $n$ variables over $R$ and is of the form $B[I^{-1}]$,
where $B$ is the symmetric algebra of a projective $R$-module $Q$ of rank $n$,
$Q$ is a direct sum of finitely generated projective $R$-modules of rank one,
and $I$ is an invertible ideal of $B$. 
\end{thm}

\section{The Epimorphism Problem}

The  Epimorphism Problem for hypersurfaces
asks the following fundamental question:

\medskip

\noindent
{\bf Question 4.}  Let  $k$ be a field and $ f \in B= k^{[n]}$ for some integer $n \ge 2$. Suppose, 
$$
B/(f)\cong k^{[n-1]}
$$
Does this imply that $B= k[f]^{[n-1]}$?, i.e., is $f$ a coordinate in $B$?

\medskip

This problem is generally known as the  \textit{the Epimorphism Problem}.
It is an open problem and is 
regarded as one of the most challenging and celebrated 
problems in the area of affine algebraic geometry (see \cite{DuG}, \cite{M} \cite{R5}, \cite{RS1} for useful surveys).

The first major breakthrough on Question 4 was achieved during 1974-75, independently, 
by S.S. Abhyankar-T. Moh and M. Suzuki (\cite{AM}, \cite{Suz}). 
They showed that Question 4 has an affirmative answer when $k$ is a field of characteristic zero and $n=2$.
Over a field of positive characteristic, explicit examples of non-rectifiable epimorphisms from $k[X,Y]$ to $k[T]$
(referred to in Section 2) and hence explicit examples of  nontrivial lines had already been demonstrated by B. Segre (\cite{Se})
in 1957 and  M. Nagata (\cite{Na}) in 1971. However over a field of characteristic zero, we have the following conjecture:

\medskip
\noindent
{\it Abhyankar-Sathaye Conjecture.}  Let  $k$ be a field of characteristic zero and $ f \in B= k^{[n]}$ for some integer $n \ge 2$. Suppose  that 
$
B/(f)\cong k^{[n-1]}.
$
Then $B= k[f]^{[n-1]}$. 

\medskip

In case $n=3$ some special cases have been solved by A. Sathaye, P. Russell and D. Wright (\cite{S1}, \cite{R2}, \cite{W} and \cite{RS}). 
In \cite{S1}, Sathaye proved the conjecture for the linear planes, i.e., polynomials $F$ of the form $aZ-b$, where $a, b \in k[X,Y]$. 
This was further extended by P. Russell over fields of any characteristic. They proved that

\begin{thm}\label{sa}
 Let $F \in k[X, Y, Z]$ be such that $F= aZ-b$, where $a (\neq 0), b \in k[X, Y]$,
and $k[X, Y, Z]/(F) = k^{[2]}$. 
Then there exist $X_0, Y_0 \in k[X, Y]$ such that $k[X, Y] = k[X_0, Y_0]$ with $a \in k[X_0]$ and
$k[X, Y, Z] = k[X_0, F]^{[1]}$. 
\end{thm}

When $k$ is an algebraically closed field of characteristic $p \ge 0$, 
D. Wright (\cite{W}) proved the conjecture for polynomials $F$ of the form $aZ^m-b$ with $a, b \in k[X,Y]$,
$m \ge 2$ and $p \nmid m$.  P. Das and A.K. Dutta showed (\cite[Theorem 4.5]{DD2}) 
that Wright's result extends to any field $k$. They proved that

\begin{thm}\label{Wr}
Let $k$ be any field with ${\rm ch}~k=p (\ge 0)$
and $F = aZ^m-b \in k[X,Y,Z]$ be such that
$a (\neq 0), b \in k[X,Y]$, $m \ge 2$ and $p \nmid m$. 
Suppose that $k[X, Y, Z]/(F) = k^{[2]}$. 
Then there exists $X_0 \in k[X, Y]$ such that $k[X, Y] = k[X_0, b]$ with $a \in k[X_0]$ and
$k[X, Y, Z] = k[F, Z, X_0]$. 
\end{thm}

The condition that $p \nmid m$ is necessary in  Theorem \ref{Wr} (cf. \cite[Remark 4.6]{DD2}).

Most of the above cases are covered by the following 
generalisation due to Russell and Sathaye (\cite[Theorem 3.6]{RS}):   
 
\begin{thm}\label{rs}
Let $k$ be a field of characteristic zero and let 
$$
F = a_mZ^m + a_{m-1} Z^{m-1} + \cdots + a_1Z + a_0 \in k[X, Y, Z]
$$
where $a_0, \dots, a_m\in k[X, Y]$ are such that ${\rm GCD}~(a_1, \dots, a_m)\notin k$. 
Suppose that $k[X, Y, Z]/(F) = k^{[2]}$. 
Then there exists $X_0 \in k[X, Y]$ such that $k[X, Y] = k[X_0, b]$ with $a_m \in k[X_0]$. 
Further,  $k[X, Y, Z]= k[F]^{[2]}$. 
\end{thm}

Thus, for $k[X,Y,Z]$, the Abhyankar Sathaye conjecture remains open for the case 
when  ${\rm GCD}~(a_1, \dots, a_m)=1$. 

A common theme in most of the partial results proved in the Abhyankar-Sathaye conjecture for $k[X,Y, Z]$ is that, 
if $F$ is considered as a polynomial in $Z$, then the coordinates of $k[X,Y]$ can be so chosen that the coefficient of 
$Z$ becomes a polynomial in $X$. The Abhyankar Sathaye conjecture for $k[X,Y,Z]$ can now be split into two parts. 

\medskip
\noindent
{\bf Question 4A.}
Let $k$ be a field of characteristic zero and let 
$$
F = a_mZ^m + a_{m-1} Z^{m-1} + \cdots + a_1Z + a_0 \in k[X, Y, Z]
$$
where $a_0, \dots, a_m\in k[X, Y]$.  
Suppose that $k[X, Y, Z]/(F) = k^{[2]}$. 
Does there exist $X_0 \in k[X, Y]$ such that $k[X, Y] = k[X_0]^{[1]}$ with $a_m \in k[X_0]$?

\medskip
\noindent
{\bf Question 4B.}
Let $k$ be a field of characteristic zero and suppose
$$
F = a_m(X)Z^m + a_{m-1} Z^{m-1} + \cdots + a_1Z + a_0 \in k[X, Y, Z]
$$
where $a_0, \dots, a_{m-1}\in k[X, Y]$ and $a_m \in k[X]$.  
Suppose that $k[X, Y, Z]/(F) = k^{[2]}$. 
Does this imply that $k[X,Y, Z]= k[F]^{[2]}$?

\medskip

L.M. Sangines Garcia in his Ph.D. thesis (\cite{San}) answered Question 4A affirmatively for the case $m=2$.
In \cite{BG4}, Bhatwadekar and the author have given  an alternative proof of this result of Garcia. 
 
\medskip

When $k$ is any field, as a partial generalisation  of Theorem \ref{sa} and Question 4B in four variables, 
the author proved the Abhyankar-Sathaye conjecture for a polynomial $F$ of the form $X^mY- F(X,Z,T)\in k[X,Y,Z,T]$.
This was one of the consequences of her general investigation on the ZCP (\cite{G2}).
In the process, she related it with other central problems on affine spaces like the affine fibration problem 
and the ZCP. The author has proved equivalence of ten statements, 
some of which involve an invariant introduced by H. Derksen which is called the Derksen invariant. 

The Derksen invariant of an integral domain $B$, denoted by $\dk(B)$, is defined as the smallest
subring of $B$ generated by the kernel of $D$, where $D$ varies over the set of all locally nilpotent derivations of $B$.  

\begin{thm}\label{ngc}
Let $k$ be a field of any characteristic and $A$ an integral domain defined by 
$$
A= k[X,Y,Z,T]/(X^m Y - F(X, Z, T)), {\text{~~where~~}} m > 1.
$$
Let $x$, $y$, $z$ and $t$ denote, 
respectively, the images of $X$, $Y$, $Z$ and $T$ in $A$. 
Set $f(Z, T):= F(0, Z, T)$ and $G:= X^m Y - F(X, Z, T)$.
Then the following statements are equivalent:
\begin{enumerate}
\item [\rm (i)] $k[X, Y, Z, T]= k[X, G]^{[2]}$. 
\item [\rm (ii)] $k[X, Y, Z, T]= k[G]^{[3]}$.
\item [\rm (iii)] $A = k[x]^{[2]}$.
\item [\rm (iv)] $A = k^{[3]}$. 
\item[\rm(v)] $A^{[\ell]}\cong_k k^{[\ell+3]}$ for some integer $\ell \ge 0$ and  $\dk(A) \neq k[x, z, t]$.
\item [\rm(vi)] $A$ is an $\A^2$-fibration over $k[x]$ and $\dk(A) \neq k[x, z, t]$.
\item[\rm(vii)] $A$ is geometrically factorial over $k$, $\dk(A)\neq k[x, z, t]$ and 
the canonical map $k^* \rightarrow K_1(A)$ (induced by the inclusion $k \hookrightarrow A$)
is an isomorphism.
\item[\rm(viii)] $A$ is geometrically factorial over $k$, $\dk(A)\neq k[x, z, t]$ and $(A/xA)^*= k^*$.
\item [\rm (ix)] $k[Z, T]= k[f]^{[1]}$. 
\item [\rm(x)] $k[Z, T]/(f)=k^{[1]}$ and $\dk(A) \neq k[x, z, t]$.
 \end{enumerate}
\end{thm}

The equivalence of (ii) and (iv) provides an answer to Question 4 for the special case of the polynomial 
$X^m Y - F(X, Z, T)$. The equivalence of (i) and (iii) provides an answer to a special case of Question 4$^\prime$ (stated below)
for the ring $R=k[x]$. 
The equivalence of (iii) and (vi) answers Question 3 in a special situation.  For more discussions, see \cite{G5}. 

In a remarkable paper S. Kaliman proved the following result over the field of complex numbers (\cite{Ka}).
Later D. Daigle and S. Kaliman extended it over any field $k$ of characteristic zero (\cite{DK}).

\begin{thm}\label{K1}
Let $k$ be a field of characteristic zero.
Let $F \in k[X, Y, Z]$ be such that $k[X, Y, Z]/(F- \lambda) = k^{[2]}$ for almost every $\lambda \in k$.
Then $k[X, Y, Z]= k[F]^{[2]}$. 
\end{thm}

A general version of Question 4 can be asked as: 

\medskip

\noindent
{\bf Question 4$^{\prime}$.}  Let  $R$ be a ring and $ f \in A= R^{[n]}$ for some integer $n \ge 2$. Suppose, 
$$
A/(f)\cong R^{[n-1]}
$$
Does this imply that $A= R[f]^{[n-1]}$?, i.e., is $f$ a coordinate in $A$?

\medskip

There have been affirmative answers to Question 4$^{\prime}$ in special cases 
by S.M. Bhatwadekar, A.K. Dutta and P. Das  (\cite{Bget}, \cite{BDl}, \cite{DD2}).
Bhatwadekar and Dutta had considered linear planes, i.e., 
polynomials $F$ of the form $aZ-b$, where $a, b \in R[X,Y]$ over a discrete valuation ring $R$ 
and proved that special cases of the linear planes are actually variables. 
Bhatwadekar-Dutta have also shown (\cite{BDi}) that a negative answer to Question 4$^{\prime}$
in the case when $n=3$ and $R$ is a discrete valuation ring containing $\bQ$ will
give a negative answer to the affine fibration problem  (Question 3(i)) for the case $n=2$ and $d=2$.
An example of a case of linear planes which remains unsolved is discussed in Section \ref{bdex}.

\section{$\A^n$-forms} 

Let $A$ be an algebra over a field $k$. 
We say that $A$ is an $\A^n$-form over $k$ if $A\otimes_k L=L^{[n]}$ for some finite algebraic extension $L$ of $k$.
Let $A$ be an $\A^n$-form over a field $k$. 

When $n=1$, it is well-known that if $L|_k$ is a separable extension, then $A=k^{[1]}$ (i.e., trivial) and 
that if $L|_k$ is purely inseparable  then $A$ need not be $k^{[1]}$.
An extensive study of such purely inseparable algebras was made by T. Asanuma in \cite{Ap}. 
Over any field of positive characteristic, the non-trivial purely inseparable $\A^1$-forms
can  be used to give examples of non-trivial $\A^n$-forms  for any integer $n>1$. 

When $n=2$ and $L|_k$ is a separable extension, then 
T. Kambayashi established  that $A=k^{[2]}$ (\cite{K}).
However,  the problem of existence of non-trivial separable $\A^3$-forms 
is open in general. A few recent partial results on the triviality of 
separable $\A^3$-forms are mentioned below.

\medskip

Let $A$ be an $\A^3$-form over a field $k$ of characteristic zero and $\bar{k}$ be an algebraic closure of $k$.
Then $A= k^{[3]}$ if it satisfies any one of the following.

(1)   $A$ admits a fixed point free locally nilpotent derivation $D$ (D. Daigle and S. Kaliman \cite[Corollary 3.3]{DK}).

(2) $A$ contains an element $f$ which is a coordinate of $A \otimes_k \bar{k}$ (Daigle and Kaliman \cite[Proposition 4.9]{DK}). 

(3) $A$  admits an effective action of a reductive algebraic $k$-group of positive dimension  (M. Koras and P. Russell \cite[Theorem C]{RK}). 

(4)  $A$ admits either a fixed point free locally nilpotent derivation or a non-confluent action of a unipotent group of dimension two (R.V. Gurjar, K. Masuda and M. Miyanishi \cite{GMM}).

(5)  $A$ admits a locally nilpotent derivation $D$ such that ${\rm rk~}(D\otimes 1_{\bar{k}})\leq{2}$ (A.K. Dutta, N. Gupta and A. Lahiri \cite{DGL}).

\medskip

Now let $R$ be a ring containing a field $k$. An $R$-algebra $A$ is said to be an $\A^n$-form over $R$ with respect to $k$ if
$A\otimes_k \bar{k}=(R\otimes_k \bar{k})^{[n]}$, where $\bar{k}$ denotes the algebraic closure of $k$.
A few results on triviality of separable $A^n$-forms over a ring $R$ are listed below. 

\medskip

Let $A$ be an $\A^n$-form over a ring $R$ containing a field $k$ of characteristic $0$. Then:

(1) If $n=1$,  then $A$ is isomorphic to the symmetric algebra of a finitely generated rank one projective module over $R$
(\cite[Theorem 7]{Du}). 

(2) If $n=2$ and $R$ is a PID containing $\bQ$, then $A= R^{[2]}$ (\cite[Remark 8]{Du}).   

(3) If $n=2$, then $A$ is an $\A^2$-fibration over $R$. 

(4) If $n=2$ and $R$ is a one-dimensional Noetherian domain,  
then there exists a finitely generated rank one projective $R$-module $Q$ such that $A \cong ({\rm Sym}_{R}(Q))^{[1]}$ (\cite[Theorem 3.7]{DGL}).   

(5) If $n=2$ and $A$ admits has a fixed point free locally nilpotent $R$-derivation over any ring $R$,
then there exists a  finitely generated rank one projective $R$-module $Q$ such that $A \cong ({\rm Sym}_{R}(Q))^{[1]}$ (\cite[Theorem 3.8]{DGL}).

\medskip

The result (3) above shows that an affirmative answer to the $\A^2$-fibration problem (Question 3 (i)) 
will ensure an affirmative answer to the problem of $\A^2$-forms over general rings. 
Over a field $F$ of any characteristic, P. Das  has shown (\cite{PD}) that any factorial $\A^1$-form $A$ over a ring $R$ 
containing $F$ is trivial if there exists a retraction map from $A$ to $R$.

We can't say much about  $\A^3$-forms over general rings till the time we solve it over fields. 

\section{An example of Bhatwadekar and Dutta}\label{bdex}

The following example arose from the study of linear planes over a discrete valuation ring by Bhatwadekar and Dutta \cite{BDi}. 
Question 5 stated below is an open problem for at least three decades. 
Let 
$$
A= \bC[T, X, Y, Z]   \text{~~and~~} R= \bC[T, F] \subset A,
$$
where $F= T X^2Z+X+T^2Y+T X Y^2$.

Let 
$$
P:= XZ+Y^2
$$
$$
G:= TY + XP
$$
and 
$$
H:= T^2Z-2TYP -XP^2
$$
Then, we can see that 
$$
XH+G^2= T^2P
$$ and $F= X+ TG$. Clearly $\bC[T,T^{-1}][F,G, H]\subseteq \bC[T, T^{-1}][X,Y,Z]$. 

Then the following statements hold.
\begin{enumerate}
\item [\rm (i)] $\bC[T, T^{-1}][X,Y,Z]=\bC[T, T^{-1}, F, G,H]=\bC[T,T^{-1}][F]^{[2]}$.
\item [\rm (ii)] $\bC[T,X,Y,Z]$ is an $\A^2$-fibration over $\bC[T, F]$.  
\item [\rm (iii)] ${\bC[T,X,Y,Z]}^{[1]}=\bC[T, F]^{[3]}$.
\item [\rm (iv)] $\bC[T,X,Y,Z]/(F) = \bC[T]^{[2]}=\bC^{[3]}$.
\item [\rm (v)] $\bC[T,X,Y,Z]/(F-f(T)) = \bC[T]^{[2]}$ for every polynomial $f(T) \in \bC[T]$.
\item[\rm(vi)] $\bC[T,X,Y,Z][1/F]=\bC[T, F,1/F, G]^{[1]}$.
\item[\rm(vii)] For any $u \in (T, F)R$, $A[1/u]= R[1/u]^{[2]}$, i.e., 
$\bC[T,X,Y,Z][1/u]=\bC[T, F,1/u]^{[2]}$.
\end{enumerate}

\medskip

\noindent
{\bf Question 5.} (a)  Is $A= \bC[T, F]^{[2]}(= R^{[2]})$?  \\
(b) At least is $A= \bC[F]^{[3]}$?

\medskip

If the answer is no to (a), then it is a counter-example to the following problems:
\begin{enumerate}
\item $\A^2$-fibration Problem over $\bC^{[2]}$ by (ii).
\item Cancellation Problem over $\bC^{[2]}$ by (iii).
\item Epimorphism problem over the ring $\bC[T]$ (see Question 4$^\prime$) by (iv).
\end{enumerate}

If the answer is no to (b) and hence to (a), then it is a counter-example also to the Epimorphism Problem for $\bC^{[4]} \twoheadrightarrow \bC^{[3]}$.

Though the above properties have been proved in several places, a proof is presented below. 
A variant of the Bhatwadekar-Dutta example was also constructed by 
S. V\'{e}n\'{e}reau in his thesis (\cite{V0}); for a discussion
on this and related examples, see  \cite{DF}, \cite{Fr} and \cite{L}.

\begin{proof}

(i) We show that 
\begin{equation}\label{generic}
\bC[T, T^{-1}][X,Y,Z]=\bC[T,T^{-1}][F,G, H].
\end{equation}
Note that
$$
X= F-TG, ~~~P= \dfrac{XH+G^2}{T^2}
$$
$$
Y= (G-XP)/T
$$
and 
$$
Z= (H+2TYP+XP^2)/T^2
$$
and hence equation (\ref{generic}) follows.

\smallskip

(ii) Clearly $A$ is a finitely generated $R$-algebra. It can be shown by standard arguments that $A$ is a flat $R$-algebra (\cite[Theorem 20.H]{Mat}).
We now show that $A \otimes_R k(\p)= k(\p)^{[2]}$ 
 for every prime ideal $\p$ of $R$. We note that $F-X \in TA$ and hence the image of $F$ in $A/TA$ is same as that of $X$. 
 Now let $\p$ be a prime ideal of $R$. Then either $T \in \p$ or $T \notin \p$. If $T \in \p$,  
 then $A \otimes_R k(\p)= k(\p)[Y,Z]= k(\p)^{[2]}$. If $T \notin \p$, then image of $T$  in $k(\p)$ is a unit 
 and the result follows from (i). 
 
(iii) Let $D= A[W]= \bC[T, X, Y, Z, W]= \bC^{[5]}$. We shall show that $D= \bC[T, F]^{[3]}= R^{[3]}$. Let 
$$
W_1:=TW+P 
$$$$
G_1:= \frac{(G-FW_1)}{T}= Y-XW-(TY+XP)(TW+P)= Y-XW-GW_1
$$
$$
H_1:= \frac{\{H+2GW_1-(F-GT)W_1^2\}}{T^2}=Z+2YW-XW^2
$$ 
Now let 
$$
G_2:= G_1+FW_1^2= (Y-XW)-TW_1(Y-XW-GW_1)=Y-XW-TW_1G_1
$$
and 
$$
W_2:= \frac{W_1-(H_1F+G_2^2)}{T}= W+2G_1W_1(Y-XW)-GH_1-TG_1^2W_1^2
$$
Then, it is easy to see that 
$$
\begin{array}{lll}
D[T^{-1}]&=& \bC[T,T^{-1}][ X, Y, Z, W]\\
&=& \bC[T,T^{-1}][F,G,H, W_1]\\
&=& \bC[T,T^{-1}][F, G_1, H_1, W_1]\\
&=& \bC[T,T^{-1}][F, G_2, H_1, W_2]
\end{array}
$$
and that $\bC[T, F, G_2, H_1, W_2] \subseteq D$. 
Let $D/TD=\bC[x,y,z,w]$, where $x, y, z, w$ denote the images of $X,Y, Z, W$ in $D/TD$.  
We now show that $D\subseteq \bC[T,F, G_2, H_1, W_2]$. 
For this it is enough to show that the kernel of  the natural map $\phi: \bC[T,F, G_2, H_1, W_2] \to D/TD$ is generated by $T$. 
We note that the image of $\phi$ is 
$$
\bC[x, y-xw, z+2yw-xw^2, w+2p(y-xw-xp^2)(y-xw)-xp(z+2yw-xw^2)]
$$ which is of transcendence degree $4$ over $\bC$. 
Hence the kernel of $\phi$ is a prime ideal of height  one and is generated by $T$. Therefore, 
$D= \bC[T,F, G_2, H_1, W_2]$.

(iv) and (v)  Let $B=\bC[T,X,Y,Z]/(F- f(T))$ for some polynomial $f \in \bC[T]$ and $S=\bC[T]$. By (ii) it follows that $B$ is an $\A^2$-fibration over $S$.
Hence, by Sathaye's theorem (\cite{S}),  $B$ is locally a polynomial ring over $S$ and hence by Theorem \ref{bcw}, $B$ is a polynomial ring over $S$.

(vi) Let $H_1:= \dfrac{FH+G^2}{T}$. Then
$$
H_1=\frac{(X+TG)(T^2Z-2TYP-XP^2)+(TY+XP)^2}{T}=TP+GH
$$
Let $H_2:= \dfrac{FH_1+G^3}{T}$. Then 
$$
\begin{array}{lll}
H_2&=&\dfrac{(X+TG)(TP+GH)+G^3}{T}\\
&= &\dfrac{T(G^2H+TGP+XP)+ G(XH+G^2)}{T}\\
&=&\dfrac{T(G^2H+TGP+XP)+ GT^2P}{T}\\
&=& G^2H+XP+ 2TGP
\end{array}
$$
Let $H_3:= \dfrac{F(H_2-G)+G^4}{T}$. Then 
$$
\begin{array}{lll}
H_3&=&\dfrac{F(G^2H+XP+2TGP-XP-TY)+G^4}{T}\\
&= &\dfrac{F(2TGP-TY)+ G^2(FH+G^2)}{T}\\
&=&\dfrac{TF(2GP-Y)+ TH_1G^2}{T}\\
&=& F(2GP-Y)+H_1G^2
\end{array}
$$
Now it is easy to see that
$$
\begin{array}{lll}
\bC[T,X,Y,Z, F^{-1}][T^{-1}]&=& \bC[T,T^{-1}][F, F^{-1}, G, H]\\
&=& \bC[T,T^{-1}][F, F^{-1}, G, H_1]\\
&=& \bC[T,T^{-1}][F, F^{-1}, G, H_2]\\
&=& \bC[T,T^{-1}][F, F^{-1}, G, H_3].
\end{array}
$$
and that the image of $\bC[T, F, F^{-1}, G, H_2]$ in $A[F^{-1}]/TA[F^{-1}]$ is of transcendence degree $3$. Hence 
$A[F^{-1}]=\bC[T, F, F^{-1}, G, H_3]=\bC[T, F, F^{-1}, G]^{[1]}$.

(vii) Let $\m$ be any maximal ideal of $R$ other than $(T, F)$. Then either $T \notin \m$ or $F \notin \m$. Thus, 
in either case, from (i) and (vi), we have  $A_\m= R_{\m}^{[2]}$. 

Let $u \in (T,F)R$. Then a maximal ideal of $R[1/u]$ is an extension of a maximal ideal of $R$ other than $(T,F)R$. 
Hence $A[1/u]$ is a locally polynomial ring in two variables over $R[1/u]$. Further any projective module over $R[1/u]$ is free. 
Thus, by Theorem \ref{bcw}, we have $A[1/u]=R[1/u]^{[2]}$. 
\end{proof}

\medskip

\noindent
{\bf Acknowledgement} The author thanks Professor Amartya Kumar Dutta for introducing and guiding her to 
this world of affine algebraic geometry. The author also thanks him for carefully going through this draft and improving the exposition.

{\small{

}}

\end{document}